\documentclass[]{article}

\usepackage{amsmath}
\usepackage{graphicx} 
\usepackage{tikz}
\usetikzlibrary{patterns}
\usepackage{comment}
\usepackage{siunitx}
\usepackage{soul}
\usepackage{dsfont}
\usepackage[authoryear]{natbib}
\usepackage{newtxtext}
\usepackage[subscriptcorrection]{newtxmath}

\usepackage{authblk} 

\graphicspath{{./art/}}

\usepackage[plain]{algorithm2e} 

\DeclareMathOperator{\lrs}{\lambda_{\mbox{\scriptsize LR}}}

\title{A note on the asymptotic distribution of the Likelihood Ratio Test statistic under boundary conditions}

\author[1]{Clara Bertinelli Salucci}
\author[1]{Anders Kvellestad}
\author[1]{Riccardo De Bin}
\affil[1]{Faculty of Mathematics and Natural Science, University of Oslo, Oslo, Norway\\
email: \texttt{clarabe@math.uio.no}}

\begin{document}
\maketitle

\section*{Abstract}
In the context of likelihood ratio testing with parameters on the boundary, we revisit two situations for which there are some discrepancies in the literature: the case of two parameters of interest on the boundary, with all other parameters in the interior, and the case where one of the two parameters on the boundary is a nuisance parameter. For the former case, we clarify that two seemingly conflicting results are consistent upon closer examination. For the latter, we clarify the source of the discrepancy and explain the different findings. As for this case the closed-form expression is valid only under positive correlation, we further propose a heuristic modification to the asymptotic distribution of the likelihood ratio test that extends its applicability to cases involving negative correlation.

\vspace{.3cm}

{
\noindent \emph{Keywords:} Chi-squared distribution mixtures; Nuisance parameters; Hypothesis testing
}

\section{Introduction}

When performing likelihood ratio tests for model comparison, we often rely on Wilks' theorem \citep{wilks}, which states that the likelihood ratio test statistic, hereafter $\lrs$, is chi-squared ($\chi^2$) distributed with degrees of freedom corresponding to the number of parameters under test. This result relies on a set of sufficient regularity conditions, including smoothness, identifiability, and the assumption that the true parameter lies in the interior of the parameter space \citep{brazzale}. Violations of the latter condition arise frequently in real-world applications, when the null hypothesis fixes a parameter at the boundary of its domain: in such cases, the maximum-likelihood estimator is constrained to lie on one side of that boundary. This asymmetry leads to nonstandard limiting behaviour for estimators and test statistics.

The seminal work on the asymptotic distribution of $\lrs$ under boundary violations dates back to \cite{chernoff}: the author showed that testing whether a parameter of interest $\theta$ lies on one side of a smooth ($p-1$)-dimensional surface in a $p$-dimensional parameter space leads to a 50:50 mixture of point mass at 0 and $\chi_1^2$, now recognized as the $\bar{\chi}(\omega,1)$ distribution \citep{kudo} with weights (0.5, 0.5). \cite{SL} extended Chernoff's work to the general composite hypothesis testing: they approximate the constrained and unconstrained parameter spaces by cones at the true boundary point and reduce the limiting law of $\lrs$ to the difference of two quadratic projections of the Gaussian score vector. \citet{Shapiro85, unifiedTheory} also developed the asymptotic distribution of the $\lrs$ under inequality constraints, deriving general $\bar{\chi}$ limits that apply to many boundary-constrained problems, including specific settings involving more than two constrained parameters of interest. More recent algebraic treatments \citep{KS, SKF} give closed-form weights for configurations that involve a larger number of nuisance parameters on the boundaries. These studies also show that na{\"i}vely applying $\chi^2$ critical values in these settings can be anti-conservative, and confirm that having nuisance parameters on the boundary may lead to non-$\bar{\chi}$ asymptotic distributions. Additional contributions include: \cite{variance} and \cite{CR}, who apply a linear mixed-effects model to assess whether a variance component is zero; \cite{andrews}, who derives the limiting distributions of constrained one-sided Wald score and likelihood-ratio tests, applying them to GARCH(1,1) and random-coefficient models; and \cite{Sen-Silvapulle}, who survey refined likelihood-based methods under inequality constraints across parametric, semiparametric, and nonparametric settings. \cite{chen2010} and \cite{chen2017} broaden these results to pseudo-likelihood ratio tests in which the nuisance component, finite- or infinite-dimensional, is replaced by a consistent estimator. Previously, \cite{susko} proposed a conditional log-likelihood ratio test and showed that, conditional on the number of parameters on the boundary, the statistic converges to a chi-squared distribution with a data-dependent number of degrees of freedom. \cite{chen2018} derive the asymptotic distribution for the modified pseudo-likelihood ratio test under this conditional scheme, and \cite{chen2024} provide the distribution under nonconvex constraints on model parameters. Along similar lines, \cite{almohamad} introduce an adaptive likelihood ratio test for general cone constraints, replacing the $\chi^2$ mixture by a single $\chi^2$ variable with data-driven degrees of freedom. Generally speaking, whenever one or more parameters lie on the boundary of the admissible region, valid inference requires employing a tailored asymptotic distribution -- often, but not necessarily, a $\chi^2$ mixture -- whose form and weights must be carefully derived.

The availability of a closed-form expression for the asymptotic distribution of $\lrs$ under constraints has considerable practical relevance: in physics, many parameters lose physical meaning unless they are constrained to be non-negative, such as expected event counts and production cross-sections in high-energy physics \citep{algeri}. Similar constraints apply to variance components in mixed-effects models in medicine \citep{variance} and in behavioural genetics \citep{gjessing}, Poisson intensities in epidemiology \citep{counting_processes}, and many other applied settings: standard tests must be adjusted to account for the boundary effects, and  simulation-based approaches such as Monte Carlo are often too computationally demanding to be applied routinely. 

In this note, we address the two specific cases treated in \cite{SL} whose results were challenged by \cite{KS}. The first (Case 7 in Self and Liang) concerns two parameters of interest on the boundary and no nuisance parameters: we show that the apparent discrepancy between the two approaches appears to stem merely from a misunderstanding. The second (Case 8) involves one parameter of interest and one nuisance parameter, both constrained: we show that the divergence arises from a different interpretation of the set of admissible parameter values under the alternative hypothesis, and we discuss the distributions arising in each approach. Finally, we propose a small heuristic modification to the approach of Kopylev and Sinha, allowing its applicability to cases with negative correlation.

\section{Theoretical Background}

In this note we consider the situation in which exactly two parameters lie on the boundary of the parameter space. Depending on the example under consideration, these will be either both parameters of interest or one parameter of interest and one nuisance parameter. Any additional parameter, if present, lies in the interior of the parameter space and can therefore be neglected in the context of this discussion, insofar as it does not affect the asymptotic distribution of $\lrs$ \citep{SL}.  Let \(\Theta\subset\mathbb{R}^p\) be the full parameter space and \(\Theta_0\subset \Theta\) the null-hypothesis subset of dimension \(r\), i.e. the set of parameter values that satisfy the null constraints, with complement \(\Theta_1=\Theta\setminus \Theta_0\).  The true parameter satisfies $ \theta_0 \in \Theta_0$. Given \(N\) independent observations, we denote by \(L(\theta)\) the likelihood and \(\ell(\theta)=\log L(\theta)\) the log-likelihood. The score function is defined as the gradient of the log-likelihood, 
$
  \nabla \ell(\theta) = {\partial \ell(\theta)}/{\partial \theta}.
$
The Fisher information matrix at \(\theta_0\) is then defined as
$
  I(\theta_0) = -\,E\bigl[\nabla^2\ell_N(\theta)\bigr]_{\theta=\theta_0}
$, and the likelihood ratio test statistic is 
\[
  \lrs \;=\;-2\log \left( \frac{\sup_{\theta\in \Theta_0}L(\theta)}{\sup_{\theta\in \Theta_1}L(\theta)} \right).
\]
  
\noindent The pioneering work by \cite{chernoff} showed that, when \(\theta_0\) lies on the boundary, one may approximate the local geometry of \(\Theta_0\) and \(\Theta_1\) by their tangent cones \(C_0\) and \(C_1\) at \(\theta_0\). The author further demonstrated that the limiting distribution of $\lrs$ can be expressed purely in terms of projections of the score vector $Z$ onto these cones: in the simplest complementary-cone setting one obtains that $\lrs$ converges in distribution to 

\begin{equation}
  \inf_{\theta \in C_1 - \theta_0}  Q(\theta \vert Z)
  \;-\;
  \inf_{\theta \in C_0 - \theta_0} Q(\theta \vert Z),
  \label{quadratic}
\end{equation}

\noindent with $Q(\theta \vert Z) =  (Z-\theta)^T I(\theta_0) (Z-\theta)$.

\cite{SL} apply Chernoff's cone-approximation to general composite hypotheses, possibly of differing dimensions and with nuisance parameters on the boundary. They first establish, under mild regularity conditions, that a maximizer \(\hat{\theta}\) of \(\ell(\theta)\) under \(H_0\) exists, that \(\hat{\theta}\) is \(\sqrt{N}\)-consistent even when \(\theta_0\) lies on the boundary, and that \(\ell(\theta)\) admits a uniform second-order Taylor expansion. \cite{geyer} observed that the Chernoff-type expansion holds only for local maximizers, and extending it to global maximizers requires imposing Clarke's regularity conditions on the parameter space \citep{clarke}. Self and Liang further perform the spectral decomposition \(I(\theta_0) = P\Lambda P^T\), where $P\in\mathbb{R}^{p\times p}$ is the orthogonal matrix with columns equal to the eigenvectors of $I(\theta_0)$ and $\Lambda$ is the diagonal eigenvalues matrix, and transform to isotropic Gaussians so that \(\tilde{Z} =\Lambda^{1/2}P^TZ\sim N_p(0,\mathbb{I}_p)\), where $\mathbb{I}_p$ is the $p$-dimensional identity matrix: thus, they reduce the asymptotic law to the difference of squared Euclidean distances from a standard normal $\tilde{Z}$ to two convex cones,

\begin{equation}
  \inf_{\theta\in \tilde{C}_0}\|\tilde{Z}-\theta\|^2
  -\inf_{\theta\in \tilde{C}_1}\|\tilde{Z}-\theta\|^2,
\label{3.2SL}
\end{equation}

\noindent with $\tilde{C}_1 = \{ \tilde{\theta} : \tilde{\theta} = \Lambda^{1/2} P^T \theta,   \forall  \theta \in C_1 - \theta_0 \}$ and $\tilde{C}_0 = \{ \tilde{\theta} : \tilde{\theta} = \Lambda^{1/2} P^T \theta, \forall  \theta \in C_0 - \theta_0 \}$.
Thanks to their geometric interpretation, \cite{SL} provided results for a wide range of cases, including cases that do not lead to $\bar{\chi}$ asymptotic distributions. However, their solution is given in implicit form and relies on a case-by-case geometric analysis. \cite{KS} instead proceeded by explicitly solving the minimization problems for the quadratic expression of Eq. (\ref{quadratic}) with $\Sigma^{-1} = I(\theta_0)$, exploiting the correlation structure to derive explicit analytical expressions for the asymptotic distribution of $\lrs$. They provided solutions mainly for cases where both the parameters of interest and the nuisance parameters lie on the boundary; interestingly, they briefly commented on Self and Liang's solution to Case 7, highlighting what they judge as an error ``likely due to a misprint''. We address this mismatch in the following section, and further discuss Case 8 in Section \ref{1001}.

\section{Two parameters of interest on the boundary} \label{section_case7}

\cite{SL} asserted that the asymptotic distribution for the case with two parameters of interest on the boundary (Case 7) is a mixture of $\chi^2_0$, $\chi^2_1$, and $\chi^2_2$ distributions with weights $(0.5 - p_{\mbox{\tiny SL}}, \, 0.5, \, p_{\mbox{\tiny SL}})$, where  

\begin{equation}
p_{\mbox{\tiny SL}} = \frac{\arccos \left( \frac{I_{12}}{\sqrt{I_{11} I_{22}}} \right)}{2 \pi}.
\label{pSL}
\end{equation}

\noindent \cite{KS} agree on the chi-squared mixture, but argue that the weights are instead $(p_{\mbox{\tiny KS}} \,, 0.5 \,, 0.5 - p_{\mbox{\tiny KS}})$, apparently swapped, where  

\begin{equation}
p_{\mbox{\tiny KS}} = \frac{\arccos \left( \rho \right)}{2 \pi}.
\label{pKS}
\end{equation}

\noindent While asserting this, they suppose Self and Liang to have inadvertently exchanged the $\chi_0^2$ and $\chi_2^2$ weights due to a misprint. However, we show here that the definitions in Equations~\eqref{pSL} and~\eqref{pKS} differ in such a way that the weights presented in the two articles are equivalent.

\begin{figure}[htbp]
      \centering 
    \includegraphics[width=.41\linewidth]{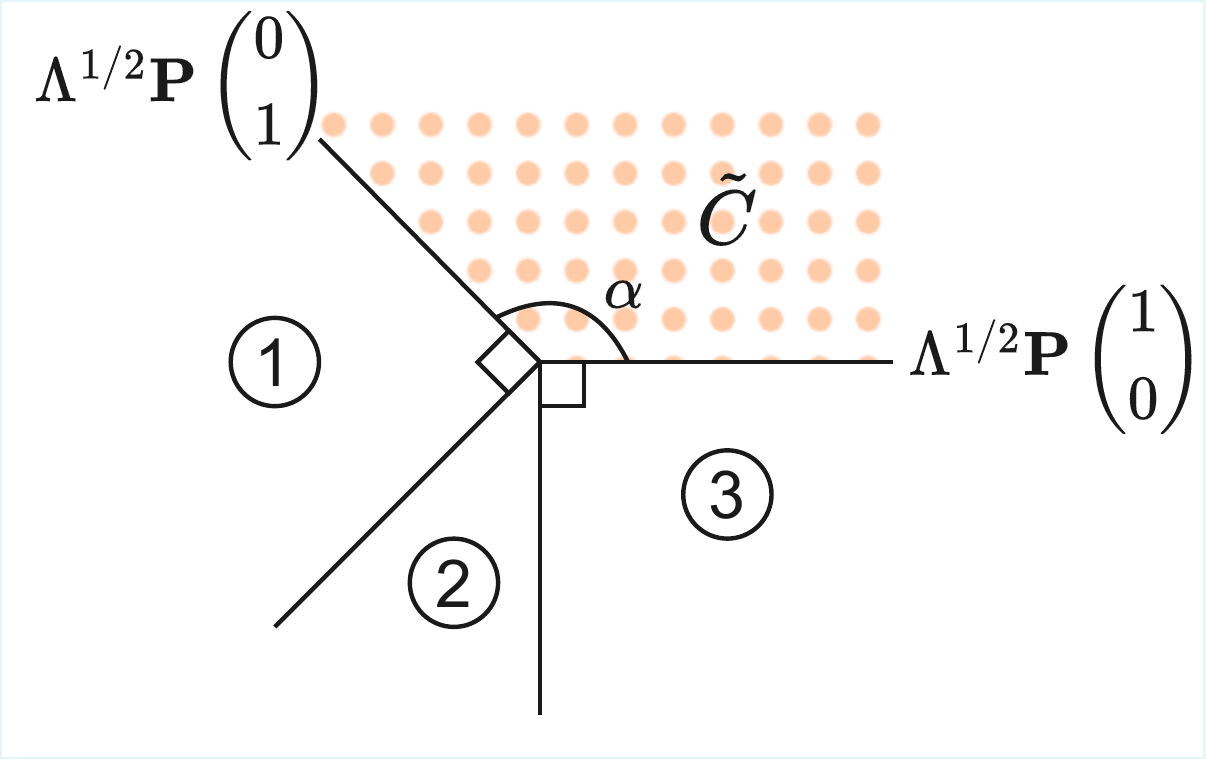}
  \caption{Adapted from \cite{SL}: diagram of the parameter space for the case of two parameters of interest on the boundary.}
  \label{SLfig1}
\end{figure}

Kopylev and Sinha seem to have supposed $\rho = {I_{12}}/{\sqrt{I_{11} I_{22}}}$, implicitly assuming that the bivariate vector $Z$ has covariance matrix $I(\theta_0)$. However, the correct covariance matrix is $\Sigma = I^{-1}(\theta_0)$ \citep[][Theorem 2]{SL}: simple algebra shows that, with the correct formulation of the covariance matrix, ${I_{12}}/{\sqrt{I_{11} I_{22}}}$ equals $-\rho$. The sign change in the arccosine argument justifies the change in the ordering of the weights: in fact, inverting the sign of the arccosine argument in Eq.~\eqref{pSL} amounts to considering another angle, say $\beta$, with  
\[
\beta = \arccos \left( -\frac{I_{12}}{\sqrt{I_{11} I_{22}}} \right) = \pi - \arccos \left( \frac{I_{12}}{\sqrt{I_{11} I_{22}}} \right) = \pi - \alpha,
\]  
which exactly maps from Region $\tilde{C}$ to Region 2 in Figure~\ref{SLfig1}. Since Region $\tilde{C}$ corresponds to the $\chi^2_2$ term, while Region 2, in which both constraints are active, to the $\chi^2_0$ term, the sign mismatch explains the mixing weights swapping and demonstrates that the two solutions are equivalent.

\section{One parameter of interest and one nuisance parameter on the boundary} \label{1001}

\subsection{Divergence in published results}

The case in which one parameter of interest and one nuisance parameter lie on the boundary is analyzed to full extent by both \cite{SL}, in their Case 8, and \cite{KS} in their Theorem 2.1. The results presented by the authors do not agree with each other: the asymptotic distribution of $\lrs$ presented by Self and Liang is a mixture of $\chi^2$ and non-$\chi^2$ distributions whereas Kopylev and Sinha, focusing only on the case with positive correlation between $Z_1$ and $Z_2$, present a simpler mixture of $\chi_0^2$, $\chi_1^2$ and $\chi_2^2$ distributions with weights  $(0.5 - q, 0.5, q)$, where $ q = \text{arcsin} \rho /(2\pi)$. 

\begin{figure}[htbp]
      \centering 
    \includegraphics[width=.49\linewidth]{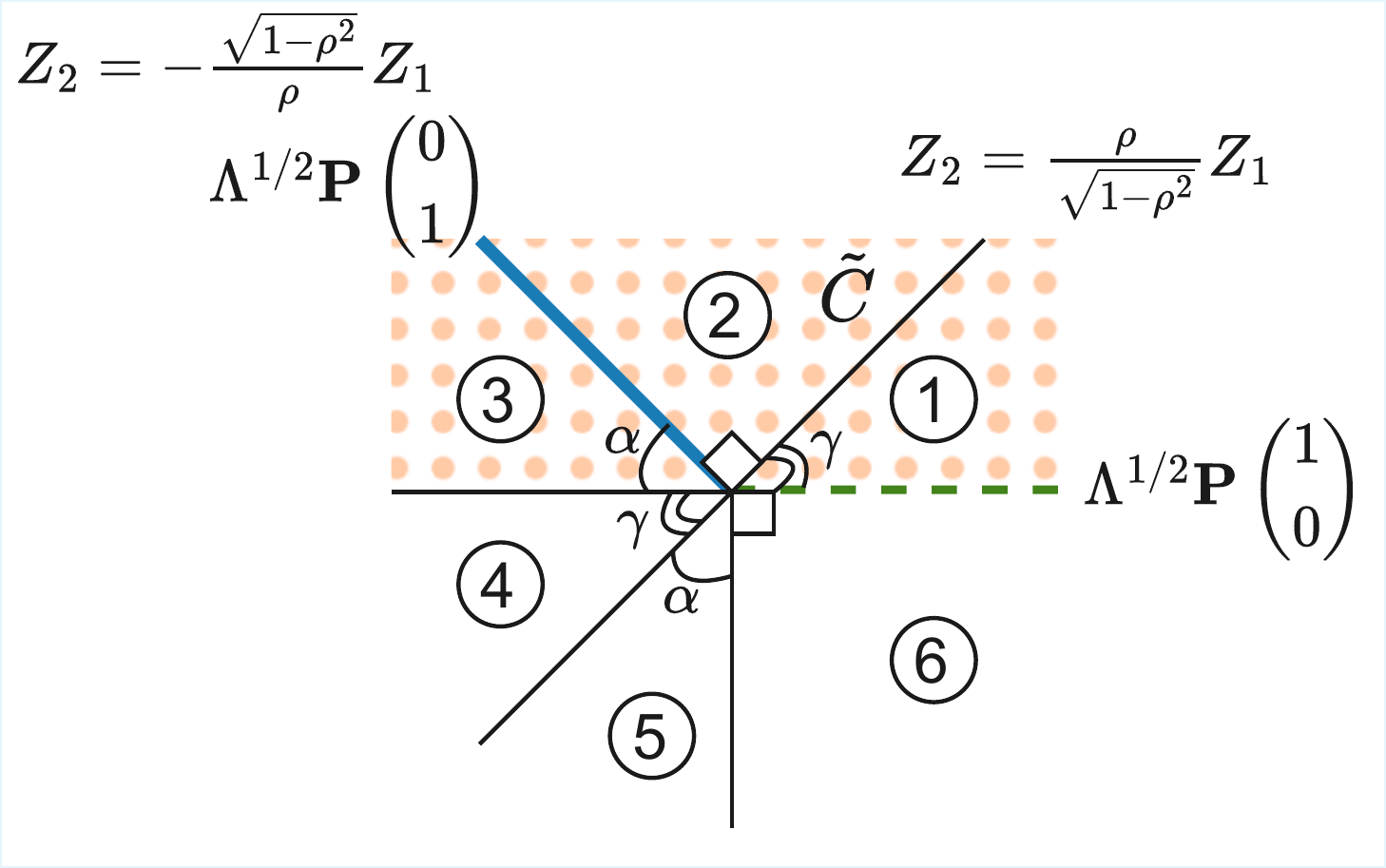}
        \includegraphics[width=.49\linewidth]{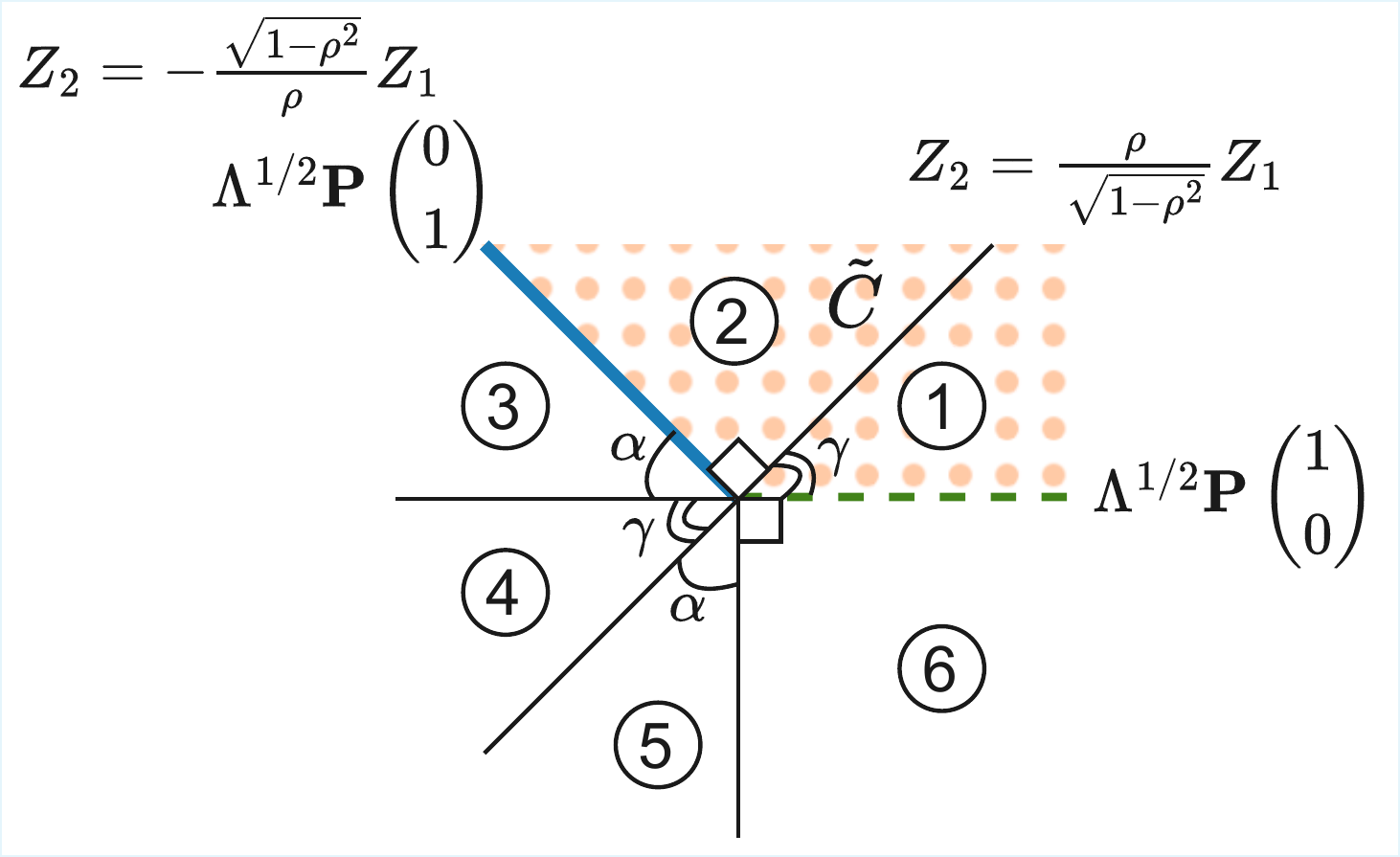}\\ 
        \hspace{1.5cm}\includegraphics[width=.4\linewidth]{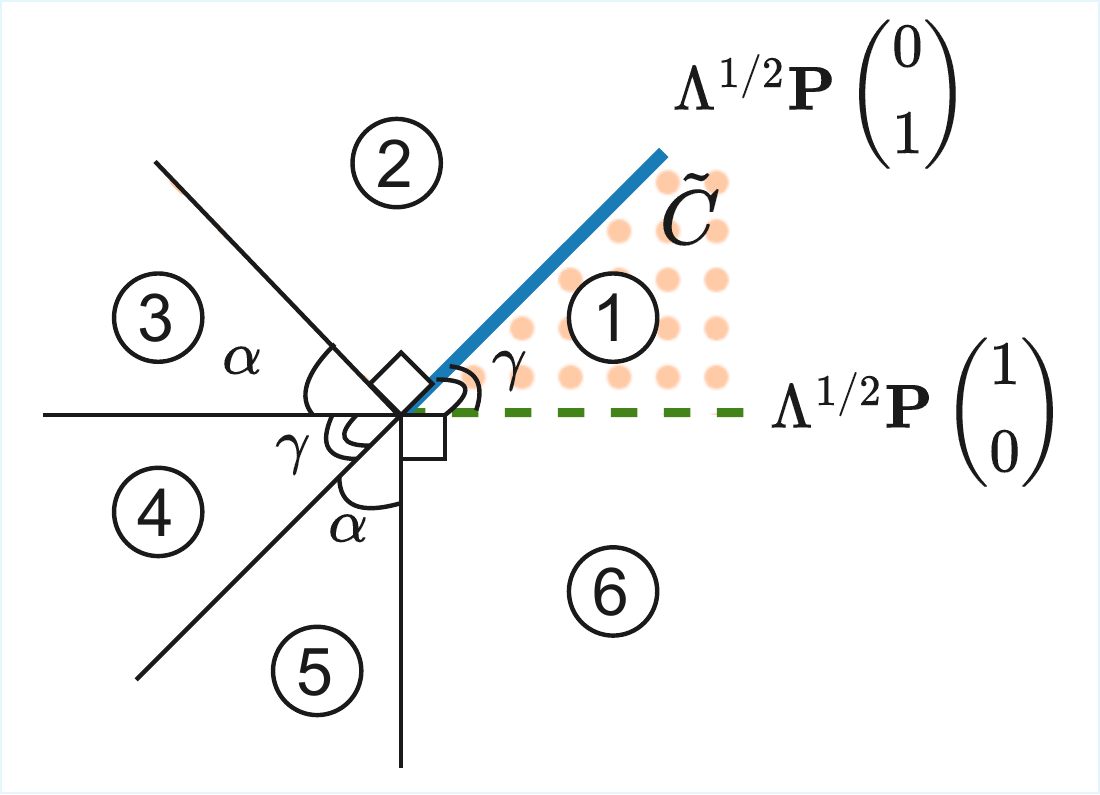}
  \caption{Parameter space for the case considered in Section 4. Top left: geometry assumed by \cite{SL}, thereby adapted. Top right: correct interpretation of the alternative cone for $\rho \geq 0$. Bottom: geometric setup induced by $\rho < 0$.}
  \label{SLfig2}
\end{figure}

The top left panel of Figure \ref{SLfig2} presents the geometry considered by Self and Liang: the axes corresponding to $\tilde{Z}_2 = 0$ and $\tilde{Z}_1 = 0$ are labelled as $\Lambda^{1/2} \mathbf{P} \left( 1 \;\; 0 \right)^\top$ (dashed) and $ \Lambda^{1/2} \mathbf{P} \left( 0 \;\; 1 \right)^\top$ (thick solid), respectively. Six regions are identified by the authors, who state that the null hypothesis cone is represented by the solid thick ray, while the alternative hypothesis cone corresponds to the dot-shaded area. Region by region, $\lrs$ is obtained from Eq. (\ref{3.2SL}) and according to the authors is: $\tilde{Z}_1^2 + \tilde{Z}_2^2$; a single squared normal random variable; $\|Y^2 \|$; $\|Y^2 \|- \tilde{Z}_1^2$; $\|Y^2 \|$; $\tilde{Z}_2^2$. Here, $Y$ is the projection of the considered point in each region onto the horizontal axis. The distribution of $\lrs$, according to these results, is chi-squared in Regions 1, 2 and 6, while non chi-squared in Regions 3, 4 and 5. Conversely, for $\rho > 0$, \cite{KS} report region-wise $\chi^2$ contributions which, taken together, yield a $\bar{\chi}$ distribution with no non-$\chi^2$ terms. 

To relate these different findings, let us denote the opening angle of Region 1 as $\gamma$, as in Figure \ref{SLfig2}. As before, we can express the $\tilde{Z}$ axes and the relative angles in terms of the untransformed score vector $Z$ and as a function of the correlation $\rho$. The angle between the $\tilde{Z}_1$- and $\tilde{Z}_2$-axes, corresponding to $\pi/2 + \gamma$, is $\arccos(I_{12}/\sqrt{I_{11} I_{22}}) = \arccos(-\rho)$. Hence, 
$
\gamma = \arccos(-\rho) - \frac{\pi}{2} = \arcsin(\rho)
$.
\noindent The slope of the thick solid line is $\tan \gamma = \rho/\sqrt(1 - \rho^2)$. Therefore, the slope of the $\tilde{Z}_2$-axis, perpendicular to it, is  $1/\tan(-\gamma) = - \sqrt(1 - \rho^2)/\rho$. Figure \ref{SLfig2} (top left panel) includes these further annotations. The complementary angle to $\gamma$, which we denote by $\alpha$, corresponds to the angle of Region 3; from the properties of complementary and supplementary angles, it follows that Regions 4 and 5 have angles $\gamma$ and $\alpha$, respectively. 

Let us assume that $\rho \geq 0$. In this case, for Region 1, the results of the two papers are in agreement: Self and Liang state that in this region $\lrs$ is $\tilde{Z}_1^2 + \tilde{Z}_2^2$, which follows a $\chi^2_2$ distribution, and the probability associated with Region 1 is $\gamma/(2\pi) = \arcsin(\rho)/(2\pi)$, which matches the result reported by Kopylev and Sinha. Similarly, both papers agree on the $\chi^2_1$ contribution which originates from Regions 2 and 6: since each of these region has a right angle, they together contribute with a weight of 1/2. The disagreement pertains Regions 3, 4 and 5. The total opening of these three regions is $\alpha + \gamma + \alpha = \pi - \gamma = \pi - \arcsin(\rho)$, which, divided by $2\pi$, yields $1/2 - \arcsin(\rho)/(2\pi)$: exactly the weight that Kopylev and Sinha associate to the $\chi_0^2$ contribution. However, Self and Liang assert that the contribution from these regions is not chi-squared, and they provide formulas that do not match any chi-squared distribution. In fact, this discrepancy arises from a misidentification of the alternative cone in \cite{SL}. The mismatch can be seen by comparing the two top panels of Figure \ref{SLfig2}, and can be understood by considering the hypotheses we are testing: under the null hypothesis, we assume that the parameter of interest is $\theta_1 = 0$ and the nuisance parameter is  $\theta_2 \geq 0$, which corresponds to the positive $\tilde{Z}_2$-axis (thick solid half-line between Regions 2 and 3), as correctly identified by Self and Liang; under the alternative hypothesis, we have both $\theta_1 \geq 0$ and $\theta_2 \geq 0$, which corresponds to the dot-shaded area in the top-right panel. Considering the whole upper half-plane as alternative cone, one includes also negative values of $\tilde{Z}_2$, violating the constraint on the nuisance parameter. When considering the correct configuration of the top-right panel, the border of the null and alternative cones to the left of the vertical $Z_2$-axis coincide, therefore the projections to the two cones from Regions 3 (now excluded from the cone), 4 and 5 coincide: one always gets $\lrs = 0$ in these regions, yielding the $\chi^2_0$ with weight $1/2 - \arcsin(\rho)/(2\pi)$ found in \cite{KS}.

Why does the chi-squared mixture not hold when $\rho < 0$? Flipping the sign of the correlation amounts to exchanging the two oblique rays in the top panels of Figure \ref{SLfig2}, so that we end up with the situation represented in the third panel of the figure: in this case, the opening angle of the alternative cone, $\gamma$, is not obtuse anymore. With $\rho > 0$, the alternative cone comprised two regions:  Region 2 is bounded by the null cone, therefore projects to it yielding a $\chi_1^2$ component, while Region 1 is far away and projects to the origin. This does not happen anymore in the new configuration ($\rho < 0$): all points in the alternative cone do have a positive projection on the null ray, therefore they contribute to a $\chi_1^2$ component, and we completely lose the $\chi^2_2$ contribution. Outside the alternative cone, all regions contribute to a $\chi_0^2$ component, because they either project to the origin or to the same point for both cones, with the exception of Region 6. A point from Region 6, in fact, projects to the thick solid ray for the null cone and the dashed horizontal ray for the alternative, yielding unavoidably a non chi-squared component. 

\subsection{A heuristic extension for negative correlation}

The ``disappearance'' of a $\chi^2_2$ region for negative correlation is corroborated by the fact that, should we blindly compute the $\bar{\chi}$ weights for $\rho < 0$, we would obtain a negative weight for the $\chi^2_2$ component. Indeed, $\arcsin(\rho) \in (-\pi/2, 0)$ for $\rho < 0$, so that $w_2 = \arcsin(\rho)/(2\pi) < 0$. As a consequence, the mixture $w_0 \, \chi_0^2 + w_1 \, \chi_1^2 + w_2 \, \chi_2^2$ with $w_0 > 0$, $w_1 > 0$ and $w_2 < 0$ is not a valid probability distribution, as it fails non-negativity. In fact, considering the case $x > 0$ (i.e. $\chi_0^2 = 0$), the mixture is

\begin{equation}
w_1 \, \chi_1^2 + w_2 \, \chi_2^2 = w_1 \left(\frac{1}{\sqrt{2\pi x}}e^{-x/2} \right) + w_2 \left(\frac{1}{2}e^{-x/2} \right) = e^{-x/2} \left( \frac{w_1}{\sqrt{2\pi x}} + \frac{w_2}{2} \right)
\label{chi_mix}
\end{equation}

\noindent which is dominated by the negative ${w_2}/{2}$ for large $x$. 

We propose an extension of the mixture cumulative distribution function (cdf) to repair the anomaly at the origin and turn the mixture into a proper probability distribution: 

\begin{equation*}
    F_{\mbox{\tiny corr}}(x) = \frac12 \, \mathds{1}(x \geq 0) + \frac12 \, F_{\chi_1^2}(x) + q \, F_{\chi_2^2}(x) - q \, G_\varepsilon(x), \qquad  \, q = \arcsin(\rho)/(2\pi) < 0 
\end{equation*}

\noindent where $G_\varepsilon(x)$ is the cdf of any non-negative density supported on $(0,\varepsilon)$ and integrating to 1.  Adding the term $-qG_\varepsilon(x)$, we restore exactly the mass that the $\chi_2^2$ term subtracted because of its negative weight. Here we take the uniform cdf on $(0, \varepsilon)$, $G_\varepsilon(x)= {x}/{\varepsilon} \;\mathds{1}(0 < x < \varepsilon)$. 
The constant $\varepsilon$ shall be chosen slightly larger than the first point where Eq. (\ref{chi_mix}) becomes non-negative.

\begin{figure}[htbp]
      \centering 
    \includegraphics[width=1\linewidth]{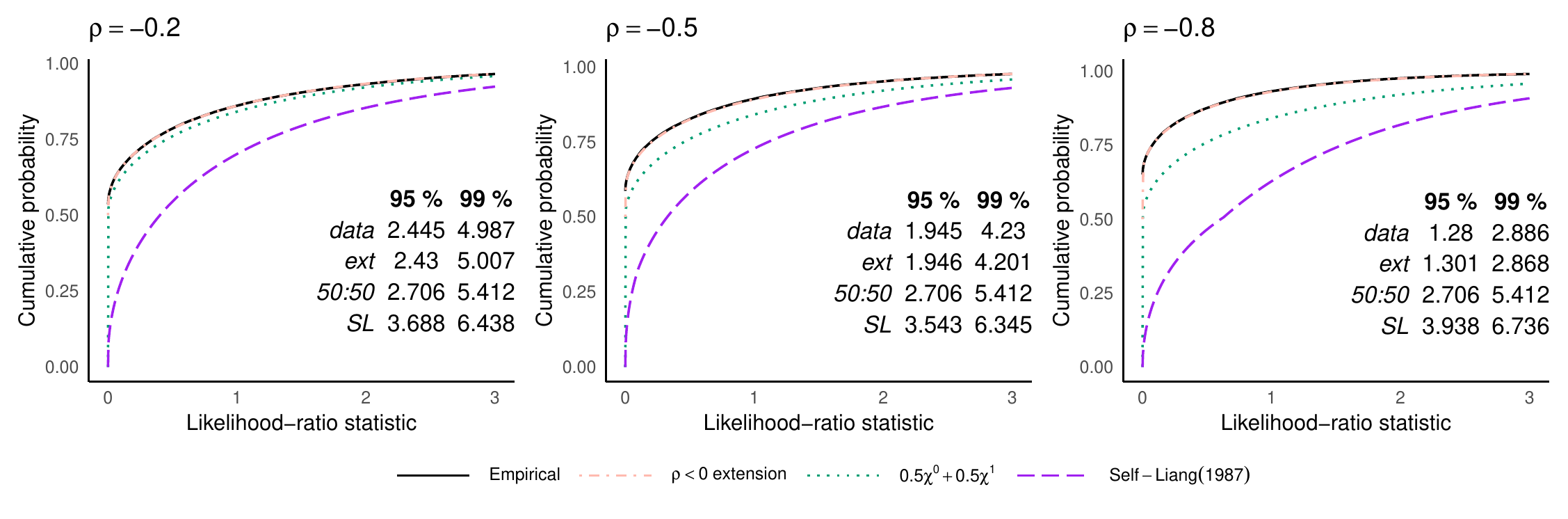}
  \caption{Empirical cdf of $\lrs$ under the null hypothesis (solid) with three analytic overlays: the proposed extension for negative correlation (dash-dotted), visually indistinguishable from the solid line; the 50:50 $\chi_0^2 + \chi_1^2$ approximation (dotted), and the solution by \cite{SL} (dashed). In-panel tables report 95\% and 99\% quantiles.}
  \label{correction}
\end{figure}

Figure \ref{correction} shows the results we achieve in a Monte Carlo (MC) study, for different values of negative correlation. The empirical cdf obtained with $10^5$ MC repetitions (sample size $N = 250$) is compared to our extension for the chi-squared mixture. Two alternative approaches are also compared: the cdf provided by \cite{SL} and the $0.5 \chi_0^2 + 0.5 \chi_1^2$ distribution, which shifts the $q$ weight from the $\chi_2^2$ component to the $\chi_0^2$ component directly. The proposed extension reproduces the empirical law with near-indistinguishable accuracy across all correlations, delivering 95\% and 99\% quantiles that coincide within MC error.
By contrast, the solution by Self and Liang sits well below the empirical cdf; the mismatch is significant at moderate correlations, and grows as the correlation becomes more negative. The alternative approach of the 50:50 mixture, instead, is quite close to the empirical distribution for small values of $\rho$, but becomes progressively conservative as $\vert \rho \vert$ increases, placing too much mass at zero and underestimating the right tail.

\section*{Acknowledgement}
This work has been supported by the Research Council of Norway (RCN) through the FRIPRO PLUMBIN' (grant n.\ 323985) and the Centre of Excellence INTEGREAT (n.\ 332645).

\bibliography{paper-ref}

\end{document}